\documentclass[reqno,12pt]{amsart}
\usepackage{latexsym}
\usepackage{amsxtra}
\usepackage{verbatim}
\usepackage{color}
\usepackage{amsthm,amsmath,amsfonts,amssymb,amscd,mathrsfs,graphics}
\usepackage{paralist}
\usepackage{amssymb}
\usepackage{times,float}
\usepackage{appendix}
\usepackage{txfonts}
\usepackage{hyperref,supertabular}
\usepackage{ifpdf}
\usepackage{enumerate}
\usepackage{fancyhdr}
\usepackage[top=1.5in, bottom=1.5in, left=1.5in, right=1.5in]{geometry}

\allowdisplaybreaks

\newcommand{\beq}{\begin{equation}}
\newcommand{\eeq}{\end{equation}}
\newcommand{\beqa}{\begin{eqnarray}}
\newcommand{\eeqa}{\end{eqnarray}}
\newcommand{\beaa}{\begin{eqnarray*}}
\newcommand{\ben}{\begin{eqnarray*}}
\newcommand{\eaa}{\end{eqnarray*}}
\newcommand{\een}{\end{eqnarray*}}

\newcommand \nc {\newcommand}

\newtheorem{theorem}{Theorem}[section]
\newtheorem{lemma}[theorem]{Lemma}
\newtheorem{proposition}[theorem]{Proposition}
\newtheorem{corollary}[theorem]{Corollary}

\nc \thref{Theorem \ref}
\nc \leref{Lemma \ref}
\nc \prref{Proposition \ref}
\nc \coref{Corollary \ref}
\nc \deref{Definition \ref}
\nc \exref{Example \ref}
\nc \reref{Remark \ref}

\newcommand{\C}{\mathbb{C}}

\newcommand{\Z}{\mathbb{Z}}

\newcommand{\f}{\mathbf{f}}

\def\d{\partial}

\def\({\left(}
\def\){\right)}
\def\[{\left[}
\def\]{\right]}
\def\<{\left\langle}
\def\>{\right\rangle}

\def\la{\lambda}

\title[The Eynard--Orantin recursion for ADE singularities]
{The Eynard--Orantin recursion for simple singularities}
\author{Todor Milanov}
\address{Kavli IPMU (WPI) \\ The University of Tokyo \\ Kashiwa \\ Chiba 277-8583 \\ Japan}
\email{todor.milanov@ipmu.jp}
\thanks{{\em 2000 Math. Subj. Class.} 14D05, 14N35, 17B69}

\begin{document}

\begin{abstract}
According to \cite{BOSS} and \cite{M1}, the ancestor correlators of
any semi-simple cohomological field theory satisfy {\em local}
Eynard--Orantin recursion. In this paper, we prove that for simple
singularities, the local recursion can be extended to a global
one. The spectral curve of the global recursion is an interesting
family of Riemann surfaces defined by the invariant polynomials of the
corresponding Weyl group. We also prove that for genus 0 and 1, the
free energies introduced in \cite{EO} coincide up to some constant
factors with respectively the genus 0 and 1 primary potentials
of the simple singularity.
\end{abstract}
\maketitle
\tableofcontents
\addtocontents{toc}{\protect\setcounter{tocdepth}{1}}

\section{Introduction}

\subsection{Motivation}
According to \cite{BOSS} and \cite{M1}, the ancestor correlators of
any semi-simple cohomological field theory satisfy {\em local}
Eynard--Orantin recursion. The term local refers to the fact that the
spectral curve is just a disjoint union of several discs. If we are
interested in computing specific ancestor Gromov--Witten (GW)
invariants in terms of Givental's $R$-matrix, then the local recursion
is all that we need. However, if we want to understand the nature of
the generating function from the point of view of representations of
vertex algebras (see \cite{BM}) and integrable systems (see
\cite{GM}), then it is important to extend the local recursion to a
global one, i.e., extend the spectral curve and the recursion kernel
to global objects (see \cite{BE}). The appropriate spectral curve however, looks quite
complicated in general, since it is parametrized by period
integrals. In particular, finding whether an appropriate generalization of
the global Eynard--Orantin recursion \cite{EO, BE} exists in the settings of
semi-simple cohomological field theories is a very challenging and
important problem. 

In this paper we would like to solve the above problem for simple
singularities. In this case, the spectral curve turns out to be a
classical Riemann surface defined by the invariant polynomials of the
monodromy group of non-maximal degree, while the invariant polynomial
of maximal degree defines a branched covering of $\mathbb{P}^1$. This
brunched covering was also studied by K. Saito (unfortunately he did
not write a text), because it is a covering of what he called a {\em primitive
direction} in the space of miniversal deformations of the
singularity.

I think that the spectral curve for simple
singularities is important also in the representation theory of the
corresponding simple Lie algebras. For example, one can obtain a
simple proof of the well known fact that the order of the Weyl group is the
product of the degrees of the invariant polynomials (see Appendix
\ref{App1}).   

Finally, after a small modification our argument should work also
for all finite reflection groups. The spectral curve is a certain
family of Hurwitz covers of $\mathbb{P}^1$ parametrized by an open
subset in the space of orbits of the corresponding reflection
group. It would be interesting to obtain the Frobenius structure on
the space of orbits of the reflection group (see \cite{Du,S2})  via the
construction of a Frobenius structure on the moduli space of Hurwitz
covers (see \cite{Du}, Lecture 5).

\subsection{Singularity theory}
Let $f\in \C[x_1,x_2,x_3]$ be a weigthed-homogeneous polynomial that
has an isolated critical point at $0$ of $ADE$ type. Such polynomials
correspond to the ADE Dynkin diagrams and are listed in Table
\ref{table1}, where we have included also the Coxeter number $h$ and
the Coxeter exponents of the corresponding simple Lie algebra. 

\begin{table}[htb]
\caption{Simple singularities}
\begin{tabular}{cllc}\label{table1}
\textbf{Type} & $\boldsymbol{f}\boldsymbol{(}\boldsymbol{x}\boldsymbol{)}$ & 
\textbf{Exponents} & $\boldsymbol{h}$ \\
\vspace{3pt}
$A_N$ & $x_0^{N+1} \!+\! x_1^2 \!+\! x_2^2$ 
& $1,2,\dots,N$  & $N\!+\!1$ \\ 
\vspace{3pt}
$D_N$ & $x_0^{N-1} \!+\! x_0 x_1^2 \!+\! x_2^2$ 
& $1,3,\dots,2N\!-\!3,N\!-\!1$ & $2N\!-\!2$ \\ 
\vspace{3pt}
$E_6$ & $x_0^4 \!+\! x_1^3 \!+\! x_2^2$ 
& $1,4,5,7,8,11$ & $12$ \\ 
\vspace{3pt}
$E_7$ & $x_0^3 x_1 \!+\! x_1^3 \!+\! x_2^2$ 
& $1,5,7,9,11,13,17$ & $18$ \\ 
\vspace{3pt}
$E_8$ & $x_0^5 \!+\! x_1^3 \!+\! x_2^2$
& $1,7,11,13,17,19,23,29$ & $30$ 
\end{tabular}
\end{table}

We fix a miniversal deformation 
\beq\label{unfolding}
F(t,x) = f(x)+\sum_{i=1}^N t_i v_i(x),\quad t=(t_1,\dots,t_N)\in
B:=\mathbb{C}^N,
\eeq
where $\{v_i(x)\}_{i=1}^N$ is a set of weighted-homogeneous
polynomials that represent a basis of the Jacobi algebra
\ben
H=\mathbb{C}[x_1,x_2,x_3]/(f_{x_1},f_{x_2},f_{x_3}).
\een
The form $\omega:=dx_1dx_2dx_3$ is primitive in the sense of K. Saito
\cite{S1, MS} and the space $B$ inherits a Frobenius structure (see
\cite{He,SaT}). For some background on Frobenius structures we refer
to \cite{Du}. The Frobenius multiplication on $T_tB$ is obtained
from the multiplication in the Jacobi algebra of $F(t,\cdot)$ via the
Kodaira-Spencer isomorphism 
\ben
T_tB\cong
C[x_1,x_2,x_3]/(F_{x_1}(t,x),\dots,F_{x_3}(t,x)),\quad \partial/\partial
t_i\mapsto \partial F/\partial t_i.
\een
While the Frobenius pairing $(\ ,\ )$ on $TB$ is the residue pairing  
\ben
(\phi_1(x),\phi_2(x))_t := \frac{1}{(2\pi\sqrt{-1})^3} \int_{\Gamma}
\frac{\phi_1(x)\phi_2(x)}{F_{x_1}(t,x)\cdots F_{x_3}(t,x)} dx_1\dots dx_N,
\een
where the cycle $\Gamma$ is a disjoint union of sufficiently small
tori around the critical points of $F$ defined by equations of the type
$|F_{x_1}|=\cdots =|F_{x_3}|=\epsilon$.  In particular, we have the following identifications:
\ben
T^*B\cong TB\cong B\times T_0B\cong B\times H,
\een
where the first isomorphism is given by the residue pairing, the
second by the Levi--Civita connection of the flat residue pairing, and the last one is the
Kodaira--Spencer isomorphism 
\beq\label{KS}
T_0B \cong H,\quad \d/\d t_i\mapsto \left. \d_{t_i}F\right|_{t=0} \ {\rm mod}\ 
 ( f_{x_1},\dots,f_{x_3}).
\eeq
Let $B_{ss}\subset B$ be the subset of semi-simple points, i.e.,
points $t\in B$ such that the critical values of $F(t,\cdot)$ form a
coordinate system in a neighborhood of $t$.  For every $t\in B_{ss}$,
using Givental's higher-genus reconstruction formalism \cite{G1,G2}, we
define  ancestor correlation functions of the following form (c.f. \cite{M1})
\beq\label{correlators}
\langle a_1 \psi_1^{k_1},\dots,a_n\psi_n^{k_n}\rangle_{g,n}(t),\quad
a_i\in H,\quad k_i\in \Z_{\geq 0} (1\leq i\leq n).
\eeq
A priory, each correlator depends analytically on $t\in B_{ss}$, but it
might have poles along the divisor $B\setminus{B_{ss} }$. According to
\cite{M2} the correlation functions \eqref{correlators} extend
analytically to the entire domain 
$B$.

\subsection{The period vectors}

Put $X=B\times \mathbb{C}^3$ and $S=B\times \mathbb{C}$. Let
$\Sigma\subset S$ be the {\em discriminant} of the map
\ben
\varphi: X\to S,\quad \varphi(t,x):=(t,F(t,x)).
\een 
Removing the singular fibers $X'=X\setminus\varphi^{-1}(\Sigma)$ we
obtain a smooth fibration $X'\to S'$, where $S'=S\setminus \Sigma$,
known as the Milnor fibration. Let us denote by
$X_{t,\lambda}=\varphi^{-1}(t,\lambda)$ the fiber over $(t,\lambda)\in
S'$. The vector spaces $H^2(X_{t,\lambda};\mathbb{C})$ and
$H_2(X_{t,\lambda},\mathbb{C})$ form the so called vanishing
cohomology and homology bundles. They are equipped with flat
Gauss--Manin connections.

We fix $(0,1)\in S$ as a reference point and denote by
$\mathfrak{h}:=H^2(X_{0,1},\mathbb{C})$. The dual space
$\mathfrak{h}^*=H_2(X_{0,1},\mathbb{C})$ is equipped with a
non-degenerate intersection pairing and we denote by $(\ |\ )$ the
negative of the intersection pairing, so that $(\alpha|\alpha)=2$ for
every vanishing cycle $\alpha$. The set $R$ of all vanishing cycles
together with the pairing $(\ |\ )$ is a root system of type
$ADE$. Moreover, according to the Picard--Lefschetz theory (see
\cite{AGV}) the image of the monodromy representation   
\beq\label{mon-rep}
\pi_1(S')\to \operatorname{GL}(\mathfrak{h}^*)
\eeq
is the Weyl group of $R$, i.e., the monodromy transformation
$s_\alpha$ along a simple loop around the discriminant corresponding
to a path along which the cycle $\alpha$ vanishes is the following
reflection
\ben
s_\alpha(x) = x- (\alpha|x)\alpha,\quad \alpha\in R,\quad x\in \mathfrak{h}^*.   
\een

Let us introduce the notation $d_x$,
where $x=(x_1,\dots,x_m)$ is a coordinate system on some manifold,
for the de Rham differential in the coordinates $x$.
This notation is especially useful when we have to apply $d_x$ to
functions that might depend on other variables as well. 
The main object of our interest are the following period integrals
\beq\label{period}
I^{(k)}_\alpha(t,\lambda) = - d_t\ (2\pi)^{-1}\, \d_\la^{k+1} \
\int_{\alpha_{t,\lambda}} d_x^{-1}\omega\ \in T_t^*B\cong H,
\eeq
where $\alpha\in \mathfrak{h}$ is a cycle from the vanishing homology,
$\alpha_{t,\lambda}\in H_2(X_{t,\lambda},\mathbb{C})$ is the parallel
transport of $\alpha$ along a reference path, 
and $d_x^{-1}\omega$ is any $2$-form $\eta\in \Omega^2_{\mathbb{C}^3}$ 
such that $d_x\eta=\omega$. The periods are multivalued analytic
functions in $(t,\la)\in B\times \mathbb{C}$  with poles along the 
discriminant $\Sigma$.

\subsection{The period isomorphism}
Let us fix a coordinate system $t=(t_1,\dots,t_N)$ on $B$ defined by a
miniversal unfolding of $f$ of the type  \eqref{unfolding}. We may
assume that $v_N(x)=1$ and denote by $t-\lambda\mathbf{1}$ the point
with coordinates $(t_1,t_2,\dots,t_N-\lambda)$. Note that
$X_{t,\lambda}=X_{t-\lambda\mathbf{1},0}$, so the period vectors have
the following translation symmetry
\beq\label{tr-sym}
I^{(k)}_\alpha(t,\lambda) = I^{(k)}_\alpha(t-\lambda\mathbf{1},0).
\eeq
Sometimes we restrict the period integrals to $\lambda=0$ and it will
be convenient to use as a reference point $-\mathbf{1}\in B$. Note
that this choice is compatible with the choice of the other reference
point $(0,1)\in B\times \mathbb{C}$ in a sense that the values of the period vectors at
these two points are identified via the translation symmetry \eqref{tr-sym}. 

Now we can state the following result that goes back to Looijenga
\cite{Lo} and Saito \cite{S1}. The monodromy covering space of 
$B':=S'\cap B$ is the covering $\widetilde{B'}$of $B'$ corresponding to the kernel of
the monodromy representation \eqref{mon-rep}. It can be constructed as
the set of equivalence classes of pairs $(t,C)$, where $t\in B'$ and $C$ is a path in
$B'$  from the reference point $-\mathbf{1}$ to $t$ and the
equivalence relation $(t_1,C_1)\sim (t_2,C_2)$ is $t_1=t_2$ and
$C_1\circ C_2^{-1}$ is in the kernel of the monodromy representation
\eqref{mon-rep}.  Note that the period integrals are by definition functions on
$\widetilde{B'}$. In particular we have a well defined {\em period
  map} 
\ben
\widetilde{\Phi}: \widetilde{B'} \to \mathfrak{h}',\quad \langle\widetilde{\Phi}(C,t),\alpha\rangle :=(I_\alpha^{(-1)}(t,0),1),
\een
where $\mathfrak{h}'$ is the complement in $\mathfrak{h}$ of the
reflection hyperplanes of the roots $R$, i.e.,
$\mathfrak{h}'=\{x\in \mathfrak{h}\ |\ \langle \alpha, x\rangle \neq
0\ 
\forall \alpha\in R\}$. The first statement is that $\widetilde{\Phi}$ is an analytic
isomorphism. In particular, there is an induced isomorphism $\Phi:
B'\to \mathfrak{h}'/W$. The 2nd statement is that $\Phi$ extends
analytically across the discriminant and the extension provides an analytic 
isomorphism $B\cong \mathfrak{h}/W:=\operatorname{Spec}
(S(\mathfrak{h}^*)^W)$. Using the isomorphism $\widetilde{\Phi}$ and the
natural projection $\widetilde{B'}\to B'$ we can think of the
coordinates $t_i$ as $W$-invariant holomorphic functions on
$\mathfrak{h}'$. The 2nd statement is equivalent to saying that each
coordinate $t_i$
extends holomorphically through the reflection mirrors, the extension
is in fact a $W$-invariant polynomial in $\mathfrak{h}$,
and the ring of all $W$-invariant polynomials is  
$S(\mathfrak{h}^*)^W=\mathbb{C}[t_1,\dots,t_N].$ We refer to \cite{Lo,
S1} for the proof of all these statements.

\subsection{The spectral curve}
Let us fix a set of simple roots $\{\alpha_i\}_{i=1}^N\subset \mathfrak{h}^*$
and denote by $x=(x_1,\dots,x_N)$ the coordinate system in
$\mathfrak{h}$ corresponding to the basis of fundamental weights
$\{\omega_i\}_{i=1}^N\subset \mathfrak{h}$, i.e., 
\ben
x=\sum_{i=1}^N x_i\omega_i,\quad x_i=\langle \alpha_i,x\rangle.
\een
As explained above $t_i\in \mathbb{C}[x_1,\dots,x_N]^W$ are invariant
polynomials and since the period mapping is weighted-homogeneous,
$t_i$ are homogeneous polynomials of certain degrees $d_i$. Let us 
assume that the degrees are in an increasing order, then  the
numbers $1=d_1-1\leq d_2-1\leq \dots d_N-1=:h-1$ are known as the
Coxeter exponents (see Table \ref{table1}). Given $s\in \mathbb{C}^{N-1}$ we define the algebraic curve
$V_s\subset \mathbb{P}^N$
\ben
t_i(X_1,\dots, X_N) = s_i X_0^{d_i},\quad 1\leq i\leq N-1.
\een
As we will see later on if $s\in B_{ss}$, then $V_s$ is
non-singular. In fact, the points $s$ for which $V_s$ has
singularities are precisely the caustic $B-B_{ss}$. I am not aware if
the family of algebraic curves $V_s,s\in B_{ss}$ has an official name
attached, but since it will be the spectral curve for the EO
recursion, we will refer to it as the {\em spectral curve of the
  singularity} or just the {\em spectral curve} when the singularity
is understood from the context. 

There is a natural projection
\beq\label{cover}
\lambda: V_s\to \mathbb{P}^1,\quad [X_0,X_1,\dots, X_N]\mapsto [X_0^h,t_N(X_1,\dots,X_N)],
\eeq
which is a branched covering of degree $|W|$, where $|A|$ denotes the
number of elements of the set $A$. The branching points are
$\lambda=u_1,\dots, u_N,\infty$, where $u_i$ are the critical values
of $F(s,x)$. By definition, the
period integral $(I^{(-1)}(s,\lambda),1)$ defines locally near a
non-branching point $\lambda\in \mathbb{P}^1$ a section of the
branched covering \eqref{cover}. It follows that the set of ramification points  
\ben
\lambda^{-1}(u_i),\quad 1\leq i\leq N
\een
is precisely the intersections of $V_s$ and the reflection mirrors
\ben
\langle \alpha, X\rangle = \sum_{i=1}^N \langle\alpha,\omega_i\rangle
X_i=0,\quad \alpha\in R_+,
\een
where $R_+$ is the set of positive roots. The remaining ramification
points are $\lambda^{-1}(\infty)$. They correspond to eigenvectors of
the Coxeter transformations with eigenvalue $\eta:=e^{2\pi\sqrt{-1}/h}$:
\ben
[X_0,X_1,\dots,X_N]\in \lambda^{-1}(\infty)
\een 
if and only if $X_0=0$ and $\sum_{i=1}^N X_i\omega_i\in \mathfrak{h}$
is an eigenvector with eigenvalue $\eta$ for a Coxeter
transformation. It is easy to see that 
the ramification index of any point in $\lambda^{-1}(u_i)$ is 2, while
the ramification index of a point in $\lambda^{-1}(\infty)$ is $h$.

\subsection{The Eynard--Orantin recursion }

We make use of the following formal series
\ben
\f^\alpha(t,\la;z) = \sum_{k\in \Z} \ I^{(k)}_\alpha(t,\la)\,
(-z)^k,\quad
\phi^\alpha(t,\lambda;z) = \sum_{k\in \Z} \
I^{(k+1)}_\alpha(t,\la)(-z)^k\, d\lambda.
\een
Note that $\phi^\alpha(t,\lambda;z) =d_\lambda \f^\alpha(t,\la;z)$.
Given $n$ cycles $\alpha_1,\dots,\alpha_n$ and a semi-simple
point $s\in B_{ss}$ we define the following $n$-point symmetric forms
\beq\label{n-point-series}
\omega_{g,n}^{\alpha_1,\dots,\alpha_n}(s;\la_1,\dots,\la_n) = 
{\Big\langle} \phi^{\alpha_1}_+(s,\la_1;\psi_1),\dots,  \phi^{\alpha_n}_+(s,\la_n;\psi_n)\Big\rangle_{g,n}(s),
\eeq
where the $+$ means truncation of the terms in the series with
negative powers of $z$. The functions \eqref{n-point-series} will be
called {\em $n$-point series} of genus $g$ or simply {\em correlator
  forms}. 
The ancestor correlators \eqref{correlators} are known to be
{\em tame} (see \cite{G3}), which by definition means that they 
vanish if $k_1+\cdots +k_n>3g-3+n$. Hence the correlator
\eqref{n-point-series} is a polynomial expression of the components of
the period vectors \eqref{period}. Thanks to the translation symmetry,
we may assume that $s_N=0$, then \eqref{n-point-series}
is a meromorphic function on the spectral curve $V_s\times
\cdots\times V_s$ with possible poles at the ramification points of
the covering \eqref{cover}. 

Let us fix $s=(s_1,\dots,s_{N-1})\in \mathbb{C}^{N-1}$ and denote by
$\gamma\in \mathfrak{h}^*$ an arbitrary cycle, s.t.,
$(\gamma|\alpha)\neq 0$ for all $\alpha\in R$. We define a set of
symmetric meromorphic differentials on $V_s^n$ 
with poles along the ramification points of $V_s$
\beq\label{weyl-cor}
\omega_{g,n}(s;p_1,\dots,p_n) := \omega_{g,n}^{\gamma,\dots,\gamma}(s;\lambda_1,\dots,\lambda_n),
\eeq 
where the RHS is defined by fixing a reference path for each
$(s,\lambda_i)\in S'$, s.t., $p_i=(I^{(-1)}(s,\lambda_i),1)$. Our main
result can be stated as follows.
\begin{theorem}\label{t1}
If $s\in B_{ss}$, then the forms $\omega_{g,n},$ $2g-2+n>0$, satisfy the
Eynard--Orantin recursion associated with the branched covering
\eqref{cover} and the meromorphic function $f_\gamma:V_s\to
\mathbb{P}^1$, $f_\gamma(x):=\langle \gamma,x\rangle$. 
\end{theorem}

The recursion will be recalled later on (see Section
\ref{sec:global-EO}). We would like however to emphasize that the
Eynard--Orantin recursion in our case differs from the standard one by
the initial condition:
\beq\label{w02}
\omega_{0,2}(x,y)=\sum_{w\in W} (\gamma|w\gamma)B(x,wy),
\eeq
where $B(x,y)$ is the Bergman kernel of $V_s$ and for $x=wy$ one has
to regularize the RHS by removing an appropriate singular term (see
Section \ref{sec:global-EO}).  Let us point out that while the Bergman
kernel depends on the choice of a Torelli marking of $V_s$, i.e., a
symplectic basis of $H_1(V_s;\mathbb{Z})$, our initial condition is
independent of the Torelli marking (see Corollary \ref{K-Bergman:cor}). This fact
could be proved also directly by using the explicit formula \eqref{hol-forms} for the
holomorphic 1-forms and some standard facts for $W$-invariant
polynomials. 
Finally, let us point out that
the set of correlators \eqref{weyl-cor} determines the set
\eqref{n-point-series}, because by definition 
\ben
\omega_{g,n}^{w_1\gamma,\dots,w_N\gamma}(s;\lambda_1,\dots,\lambda_N) = 
\omega_{g,n}(s;w_1^{-1}p_1,\dots,w_N^{-1}p_N),
\een
where $w_1,\dots,w_N\in W$ are arbitrary. 

The branched covering \eqref{cover} and the meromorphic function
$f_\gamma$ determine a birational model of $V_s$ in $\mathbb{C}^2$.
Following \cite{EO} we can introduce the tau-function
$Z(\hbar,s):=Z_\hbar(V_s,\lambda,f_\gamma)$ of the birational model of
the spectral curve. It has the form
\ben
Z(\hbar,s)=\exp \Big( \sum_{g=0}^\infty \hbar^{g-1} \underline{F}^{(g)}(s)\Big),
\een
where $\underline{F}^{(g)}(s)$ is called the genus-$g$ free energy of $V_s$. 
It is very natural to compare $\underline{F}^{(g)}(s)$ with Givental's
primary genus-$g$ potentials $F^{(g)}(s)$. Unfortunately we could not solve this
problem in general, but only for $g=0$ and $g=1$
\ben
\underline{F}^{(0)}(s) = -(\gamma|\gamma)\frac{|W|}{N}F^{(0)}(s),\quad
\underline{F}^{(1)}(s) = \frac{|W|}{2} F^{(1)}(s)=0,
\een
where the first identity is valid up to quadratic terms in the
Frobenius flat coordinates  of $s\in B$, while
the second one up to a constant independent of $s$. It is known that the
genus-1 potential of the Frobenius structure is homogeneous of degree
0, so it must vanish in the case of a simple singularity.
According to \cite{EO}, $Z(\hbar,s)$ satisfies Hirota bilinear
equations. According to \cite{GM}, the total ancestor
potential $\mathcal{A}_s(\hbar;\mathbf{q})$ of an ADE singularity
satisfies the Hirota bilinear 
equations of the corresponding generalized KdV hierarchy. It will be
interesting to clarify the relation between $\underline{F}^{(g)}(s) $ and
$F^{(g)}(s)$ for $g\geq 2$, as well as to determine whether Givental's
{\em primary} ancestor potential $\mathcal{A}_s(\hbar;0)$ also
satisfies Hirota bi-linear  equations.

\section{Analytic extension of the kernel of the local Eynard--Orantin recursion}

It was proved in \cite{M1} (see also \cite{BOSS}) that the correlator
forms \eqref{n-point-series} satisfy a local Eynard--Orantin (EO)
recursion, whose kernel is defined by the symplectic pairing of
certain period series. In this section, we will prove that these
symplectic pairings are convergent and can be extended to the entire
spectral curve $V_s$. Moreover, the corresponding extensions can be expressed
in an elegant way via the so called {\em Bergman kernel} of $V_s$.  

\subsection{The kernel of the local recursion}
Recall the symplectic pairing 
\ben
\Omega(f(z),g(z))=\operatorname{Res}_{z=0}(f(-z),g(z))dz,\quad f,g\in H(\!(z^{-1})\!).
\een
The local recursion is defined in terms of the symplectic pairings 
\ben
\Omega(\mathbf{\phi}^\alpha_+(s,\lambda;z),\mathbf{f}^\beta_-(s,\mu;z))=d\lambda\,
\sum_{k=0}^\infty
(-1)^{k+1}(I^{(k+1)}_\alpha(s,\lambda),I^{(-k-1)}_\beta(s,\mu)),
\een
where the infinite series is interpreted formally in a neighborhood of
a point $\mu=u_i(s)$, s.t., the cycle $\beta_{s,\mu}$ vanishes over
$\mu=u_i(s)$. We are going to prove that this infinite series
expansion is convergent to a meromorphic function on $V_s\times V_s$.

To begin with, let us recall that the periods satisfy the following
system of differential equations
\begin{align}
\notag
\partial_a I^{(n)}(s,\lambda) & = -v_a\bullet_s I^{(n+1)}(s,\lambda) \\
\notag
\partial_\lambda I^{(n)}(s,\lambda)  & = I^{(n+1)}(s,\lambda) \\
\label{de-prim}
(\lambda-E\bullet_s)\partial_\lambda I^{(n)}(s,\lambda)   & = 
\Big(\theta -n-\frac{1}{2}\Big) I^{(n)}(s,\lambda),
\end{align} 
where $\partial_a:=\partial/\partial s_a$, $E$ is the Euler vector field $E=\sum_{i=1}^N
\operatorname{deg}(t_i)t_i\partial_{t_i}$, and $\theta$ is the
Hodge-grading operator
\ben
\theta:H\to H,\quad v_a\mapsto (D/2-\operatorname{deg}(v_a))v_a,
\een
where $D=\operatorname{deg}({\rm Hess}(f)) = 1-2/h$ is the conformal
dimension of the Frobenius structure.
The key to proving the convergence is the so called {\em phase 1-form} (see \cite{G3, BM})
\begin{eqnarray*}
\mathcal{W}_{\alpha,\beta}(s,\xi) = I^{(0)}_\alpha(s,\xi)\bullet
I^{(0)}_\beta (s,0)\quad \in \quad T_s^*S,
\end{eqnarray*}
where the period vectors are interpreted as elements in $T^*_sS$ and
the multiplication in $T^*_sS$ is induced by the Frobenius
multiplication via the natural identification $T^*_sS\cong T_sS$. The
dependence on the parameter $\xi$ is in the sense of a germ at
$\xi=0$, i.e., Taylor's series expansion
about $\xi=0$. The phase form is a power series in $\xi$ whose
coefficients are multivalued 1-forms on $B'$. 
\begin{lemma}\label{Saito-formula}
We have
\begin{eqnarray*}
(\alpha|\beta) = -\iota_{E} \mathcal{W}_{\alpha,\beta}(s,0) =
-(I^{(0)}_\alpha(s,0),E\bullet I^{(0)}_\beta(s,0)).
\end{eqnarray*}
\end{lemma}
This is a well known fact due originally to K. Saito \cite{S1}. 

\begin{lemma}\label{pf:euler-equation}
The phase form is weighted-homogeneous of weight $0$, i.e., 
\begin{eqnarray*}
(\xi\partial_\xi + L_E) \mathcal{W}_{\alpha,\beta}(s,\xi) = 0,
\end{eqnarray*} 
where $L_E$ is the Lie derivative with respect to the vector field $E$.
\end{lemma}
\proof
Note that 
\begin{eqnarray*}
\mathcal{W}_{\alpha,\beta}(s,\xi) = (I^{(0)}_\alpha(s,\xi),d I^{(-1)}_\beta(s,0)).
\end{eqnarray*}
It is easy to check that $\mathcal{W}_{\alpha,\beta} $ is a closed
1-form, so using the Cartan's magic formula $L_E=d_s\iota_E + \iota_Ed_s$,
where $\iota_E$ is the contraction by the vector field $E$, we get 
\begin{eqnarray*}
L_E \mathcal{W}_{\alpha,\beta} = d_s(I^{(0)}_\alpha(s,\xi), (\theta+1/2)
I^{(-1)}_\beta(s,0)) = -d_s((\theta-1/2) I^{(0)}_\alpha(s,\xi), 
I^{(-1)}_\beta(s,0)).
\end{eqnarray*}
We used that $\theta$ is skew-symmetric with respect to the residue pairing
and that 
\begin{eqnarray*}
\iota_E d_s I^{(-1)}_\beta(s,0) = E I^{(-1)}_\beta(s,0)) = (\theta+1/2) I^{(-1)}_\beta(s,0),
\end{eqnarray*}
where the last equality comes from the differential equation
(\ref{de-prim}) with $n=-1$ and $\lambda=0$. Furthermore, using the
Leibnitz rule we get
\begin{eqnarray*}
-((\theta-1/2) d_s I^{(0)}_\alpha(s,\xi), I^{(-1)}_\beta(s,0)) 
-((\theta-1/2) I^{(0)}_\alpha(s,\xi), d_s I^{(-1)}_\beta(s,0)).
\end{eqnarray*}
The first residue pairing is
\begin{equation}\label{res-1}
(A I^{(1)}_\alpha(s,\xi),(\theta+1/2)  I^{(-1)}_\beta(s,0))  = 
-(A I^{(1)}_\alpha(s,\xi),E\bullet I^{(0)}_\beta(s,0)),
\end{equation}
where we used that $\theta$ is skew-symmetric and that
$d_sI^{(0)}_\alpha = - AI^{(1)}_\alpha$ with $A=\sum_{a=1}^N
(\partial/\partial s_a \bullet)ds_a$. 
Similarly, the 2nd residue pairing becomes
\begin{equation}\label{res-2}
((\xi\partial_\xi + E)I^{(0)}_\alpha(s,\xi), d_sI^{(-1)}_\beta(s,0)) =
\xi\partial_\xi\mathcal{W}_{\alpha,\beta}(s,\xi)
+(E\bullet I^{(1)}_\alpha(s,\xi),A I^{(0)}_\beta(s,0)).
\end{equation}
On the other hand, since the Frobenius multiplication is commutative, 
$[A,E\bullet]=0$, so the terms (\ref{res-1}) and (\ref{res-2}) add up
to $ \xi\partial_\xi\mathcal{W}_{\alpha,\beta}(s,\xi),$ which
completes the proof. 
\qed

\medskip
Let us define the following meromorphic 1-forms on $V_s\times V_s$. Given $(x,y)\in V_s\times V_s$, s.t., $x,y$ are not ramification
points, there are unique pairs $(C,\lambda)$ and $(C',\mu)$, s.t., 
$x=\widetilde{\Phi}(C,s-\lambda\mathbf{1})$ and
$y=\widetilde{\Phi}(C',s-\mu\mathbf{1})$, where $C$ and $C'$ are paths
in $B'$ connecting respectively $s-\lambda\mathbf{1}$ and
$s-\mu\mathbf{1}$ with the reference point. Put
\ben
K_{\alpha,\beta}(s,x,y):=\frac{d\lambda}{\lambda-\mu}\,
 (I^{(0)}_\alpha(s,\lambda), (\theta+1/2) I^{(-1)}_\beta(s,\mu)),
\een
where the branches of $I^{(0)}_\alpha$ and $I^{(-1)}_\beta$ are
determined respectively by the paths $C$ and $C'$. This definition
extends analytically across the ramification points and the form has a pole
only along the divisor
\ben
\{(x,y)\in V_s\times V_s\ :\ \lambda(x)=\lambda(y)\} = \cup_{w\in W} \{x=wy\},
\een
where the action of $W$ on $\mathfrak{h}$ is induced from the action
on $\mathfrak{h}^*$, i.e., 
\ben
\langle \alpha,wx\rangle=\langle w^{-1}\alpha,x\rangle,\quad \alpha\in
\mathfrak{h}^*,\quad x\in \mathfrak{h}.
\een 
\begin{proposition}\label{phase-form}
The symplectic pairing 
\ben
\Omega(\mathbf{\phi}^\alpha_+(s,\lambda;z),\mathbf{f}^\beta_-(s,\mu;z))=
K_{\alpha,\beta}(s,x,y),
\een
where $x=\widetilde{\Phi}(C,s-\lambda\mathbf{1})$,
$y=\widetilde{\Phi}(C,s-\mu\mathbf{1})$ and $C$ is the path that
specifies the value of the symplectic pairing. 
\end{proposition}
\proof 
Using the differential equations for the periods, it is easy to verify
that 
\ben
d_s\Omega(\mathbf{\phi}^\alpha_+(s,\lambda;z),\mathbf{f}^\beta_-(s,\mu;z))
= d\lambda I^{(1)}_\alpha(s,\lambda)\bullet_s I^{(0)}_\beta(s,\mu) = 
d_\lambda \mathcal{W}_{\alpha,\beta}(s-\mu\mathbf{1},\lambda-\mu).
\een
According to Lemma \ref{pf:euler-equation} we have 
\begin{eqnarray}\label{identity-1}
\partial_\lambda \mathcal{W}_{\alpha,\beta}(s',\lambda-\mu) = 
-d_{s'}\ \Big( \frac{1}{\lambda-\mu}\, \iota_E
\mathcal{W}_{\alpha,\beta}(s',\lambda-\mu)\Big),
\end{eqnarray}
which by definition is
\begin{eqnarray*}
d_{s'}\ \Big(\frac{1}{\lambda-\mu}\,
 (I^{(0)}_\alpha(s',\lambda-\mu), (\theta+1/2) I^{(-1)}_\beta(s',0)\Big).
\end{eqnarray*}
Integrating \eqref{identity-1} with respect to $s'$ along a short path
from $s_0:=s-u_i(s)\mathbf{1}$ to $s-\mu\mathbf{1}$ and
using that $I^{(-1)}_\beta(s',0)$ vanishes as $s'\to s_0$, we get  
\begin{eqnarray}\label{kernel-local}
\Omega(\mathbf{\phi}^\alpha_+(s,\lambda;z),\mathbf{f}^\beta_-(s,\mu;z))
= 
\frac{d\lambda}{\lambda-\mu}\,
 (I^{(0)}_\alpha(s,\lambda), (\theta+1/2) I^{(-1)}_\beta(s,\mu)).
\qed
\end{eqnarray}

\subsection{The local kernel and the Bergman kernel}

Now we are in a position to prove the key result in this paper. Let us
fix a symplectic basis $\{\mathcal{A}_i,\mathcal{B}_i\}_{i=1}^g$ of
$H_1(V_s;\mathbb{Z})$, s.t., $\mathcal{A}_i\circ
\mathcal{B_j}=\delta_{i,j}$. There is a unique symmetric
differential $B(x,y)\in
\Omega_{V_s}^1\boxtimes\Omega_{V_s}^1(2\Delta)$ which is holomorphic
on $V_s\times V_s$ except for a pole of order 2 with no residue along
the diagonal $\Delta\subset V_s\times V_s$, normalized by 
\ben
\oint_{y\in \mathcal{A}_i} B(x,y) = 0,\quad 1\leq i\leq N
\een 
and 
\ben
B(x,y) = \frac{d\lambda(x)d\lambda(y)}{(\lambda(x)-\lambda(y))^2}+\cdots
\een
for any local coordinate $\lambda:U\to \mathbb{C}$ ($U\subset V_s$)
and for all $x,y\in U\times U$. The differential $B(x,y)$ is called
the {\em Bergman kernel}. We refer to \cite{EO} for more details and
references. 
\begin{proposition}\label{K-Bergman}
The following identity holds
\ben
d_y K_{\alpha,\beta}(s,x,y) = \sum_{w\in W} (\alpha|w\beta)B(x,w y),
\quad 
\forall \alpha,\beta\in \mathfrak{h}^*.
\een
\end{proposition}
\proof
Put 
\ben
x=\widetilde{\Phi}(C,s-\lambda\mathbf{1}),\quad y=\widetilde{\Phi}(C',s-\mu\mathbf{1}).
\een
By definition $wy = \widetilde{\Phi}(C'\circ
w^{-1},s-\mu\mathbf{1})$. Therefore,
if $x$ and $wy$ are near by, then $w$ must be the monodromy along
the loop $C^{-1}\circ C'$. 
Using Saito's formula (\ref{Saito-formula}) we get that the leading order term of
$K_{\alpha,\beta}(s,x,y)$ near the $w$-diagonal $x=wy$ is 
\begin{eqnarray*}
\frac{d\lambda}{(\lambda-\mu)}\,
 (I^{(0)}_\alpha(s,\mu), (\mu-E\bullet) I^{(0)}_{w\beta}(s,\mu)) =
(\alpha|w\beta)
\frac{d\lambda}{(\lambda-\mu)},
\end{eqnarray*}
where we used that 
\ben
I^{(0)}_{\beta_{C'}}(s,\mu) = I^{(0)}_{(w\beta)_C}(s,\mu),
\een
where the index in the cycle denotes the path along which the cycle
has to be transported in order to define the period.

We get that the difference of the two sides of the identity that we
want prove is a holomorphic symmetric 2-form $D(x,y)$ on $V_s\times
V_s$. To prove that such a form vanishes it is enough to prove that 
\ben
\oint_{x\in \mathcal{A}_i}D(x,y)=0,\quad \forall i=1,2,\dots,g.
\een
This is true for the Bergman kernel by definition, while for $d_y
K_{\alpha,\beta}(s,x,y)$, since it is an exact form, the corresponding integral
vanishes for all cycles $\mathcal{A}\in H_1(V_s,\mathbb{Z})$ not only
$\mathcal{A}_i$. 
\qed
\begin{corollary}\label{K-Bergman:cor}
The 2-form \eqref{w02} is independent of the choice of Torelli marking.
\end{corollary}

\section{From local to global}

In this section we prove Theorem \ref{t1}.

\subsection{The unstable range}
By definition, the ancestor potential does not have non-zero correlators in the unstable range
$(g,n)=(0,0),(0,1),(0,2)$ and $(1,0)$. However, in order to formulate the
EO recursion, it is convenient to extend the definition of the
correlators in the unstable range as well in the  following two cases:
\begin{align}\label{unst-1}
\omega_{0,2}^{\alpha_1,\alpha_2}(s;\lambda_1,\lambda_2) & 
:= \Omega(\mathbf{\phi}^{\alpha_1}_+(s,\lambda_1;z), \phi^{\alpha_2}_+(s,\lambda_2;z)_-),\\
\label{unst-2}
\omega_{0,2}^{\alpha_1,\alpha_2}(s;\lambda,\lambda)
& :=  P^{(0)}_{\alpha_1,\alpha_2}(s,\lambda),
\end{align}
where $P^{(0)}_{\alpha_1,\alpha_2}(s,\lambda)$ is defined as the limit
$\mu\to \lambda$ of 
\ben
\Omega(\mathbf{\phi}^{\alpha_1}_+(s,\lambda;z),
\phi^{\alpha_2}_+(s,\mu;z)_-) -(\alpha_1|\alpha_2)\frac{d\lambda d\mu}{(\lambda-\mu)^2}.
\een
The limit exists, because the above difference is analytic near
$\mu=\lambda$ (see \cite{M1} for more details). Let us point out that
in the definition of the correlator form
\eqref{unst-2} we assume that there is a fixed path $C$ from the
reference point to $s-\lambda\mathbf{1}$. It is more natural however
to assume that there are two such paths $C_1$ and $C_2$: one for the
1st and one for the 2nd slot of the correlator form. Since 
\ben
P_{C_2}\alpha_2 = P_{C_1} (P_{C_1^{-1}\circ C_2})\alpha_2,
\een
where $P_{C_i}:H_2(X_{-\mathbf{1},0};\mathbb{C})\to
H_2(X_{s-\lambda\mathbf{1},0};\mathbb{C})$ is the parallel transport
with respect to the Gauss--Manin connection, we get that if we want to
allow two different paths in the definition \eqref{unst-2} and still
have compatibility with the monodromy representation, then we should
define
\ben
\omega_{0,2}^{\alpha_1,\alpha_2}(s;\lambda,\lambda)
 :=  P^{(0)}_{\alpha_1,w\alpha_2}(s,\lambda),
\een
where $w=P_{C_1^{-1}\circ C_2}$ and the branch on the RHS is
determined by $C_1$.

\subsection{The local EO recursion}\label{sec:loc}
According to \cite{M1}, the ancestor correlators satisfy the following recursion
\begin{align}\notag
&&
\omega_{g,n+1}^{\alpha_0,\alpha_1,\dots,\alpha_n}(s;\lambda_0,\lambda_1,\dots,\lambda_n)=
-\frac{1}{4}\,\sum_{j=1}^N\, {\rm Res}_{\lambda=u_j}\,
\frac{ \Omega(\mathbf{\phi}^{\alpha_0}_+(s,\lambda_0;z) \, ,{\bf f}^{\beta_j}_-(s,\lambda;z)) }
 { (I^{(-1)}_{\beta_j}(s,\lambda),\mathbf{1}) d\lambda} \times \\
\label{local-EO-2}
&&
\left(
\omega_{g-1,n+2}^{\beta_j,-\beta_j,\alpha_1,\dots,\alpha_n}(s;\lambda,\lambda,\lambda_1,\dots,\lambda_n)+
\sum_{\substack{ g'+g''=g\\ I'\subset \{1,\dots,n\} }} 
\omega_{g',n'+1}^{\beta_j,\alpha_{I'}}(s;\lambda,\lambda_{I'})\,
\omega_{g'',n''+1}^{-\beta_j,\alpha_{I''}}(s;\lambda,\lambda_{I''})
\right),
\end{align}
for all stable pairs $(g,n+1)$, i.e., $2g-2+n\geq 0$, where the
notation is as follows. All unstable
correlators on the RHS are set to $0$, except for the ones of the type
(\ref{unst-1}) and (\ref{unst-2}). The summation is over all
subsets $I'\subseteq \{1,2,\dots,n\}$ and for each subset
$I'=\{i_1<\cdots<i_{n'}\}$ we put
\ben
I''=\{1,2,\dots,n\}-I'=:\{j_1<\cdots < j_{n''}\}.
\een
In particular, $ n'=|I'|$ and $ n''=|I''|$. If $x=(x_1,\dots, x_n)$ is
a sequence of $n$ elements, then we define 
\ben
x_{I'}=(x_{i_1},\dots,x_{i_{n'}}),\quad x_{I''}=(x_{j_1},\dots,x_{j_{n''}}).
\een
Finally, $\beta_j$ $(1\leq j\leq N)$ is a vanishing cycle vanishing
over $\lambda=u_j$. 

\subsection{The global EO  recursion}\label{sec:global-EO}

Let us write down the recursion from Theorem \ref{t1}. Let us denote
by 
\ben
\{ y_{j,a}\ :\ 1\leq a\leq |W|/2 \ \}:=\lambda^{-1}(u_j),\quad 1\leq j\leq N ,
\een
the ramification points on $V_s$ with ramification index 2. There is a
unique root $\beta_{j,a}\in R_+$, s.t., $\langle
\beta_{j,a},y_{j,a}\rangle =0$ and the reflection $s_{\beta_{j,a}}$
induces a deck transformation $\theta_{j,a}: V_s\to V_s$ which is a
generator for the Galois group of a neighborhood of $y_{j,a}$ viewed
as a 2-sheeted covering of a neighborhood of $u_j$. 
\begin{align}
&&
\notag
\omega_{g,n+1}(s;x_0,x_1,\dots,x_n) = \sum_{j=1}^N\sum_{a=1}^{|W|/2}
\operatorname{Res}_{y=y_{j,a}}
\frac{\frac{1}{2}\int_{y}^{\theta_{j,a}(y)}
  B(x_0,y)}{(f_\gamma(y)-f_\gamma(\theta_{j,a}(y)))d\lambda(y)} \\
\notag
&&
\Big(\omega_{g-1,n+1}(s;y,\theta_{j,a}(y),x_1,\dots,x_n) + \sum_{\substack{
    g'+g''=g\\ I'\subset \{1,\dots,n\} }} 
\omega_{g',n'+1}(s;y,x_{I'})\,
\omega_{g'',n''+1}(s;\theta_{j,a}(y),x_{I''})\Big),
\end{align}
where the summation is the same as in the local recursion (see Section
\ref{sec:loc}). Let us also point out that in the above recursion all
unstable correlators are set to 0, except for 
\ben
\omega_{0,2}(x_1,x_2) =
\omega_{0,2}^{\gamma,\gamma}(s;\lambda_1,\lambda_2),\quad x_1,x_2\in V_s,
\een
where in order to define the RHS we choose $\lambda_i=\lambda(x_i)$
and paths $C_i$ in $B'$ from $-\mathbf{1}$ to $s-\lambda_i\mathbf{1}$,
s.t., $\widetilde{\Phi}(C_i,s-\lambda_i\mathbf{1})=x_i $. Using
Proposition \ref{K-Bergman} we can express the 
form $\omega_{0,2}(x,y)$  in terms of the Bergman kernel. Namely, if
$\lambda(x)\neq \lambda(y)$, then $\omega_{0,2}(x,y)$ is given by
formula \eqref{w02}. If $\lambda(x)=\lambda(y)$, then $x=w_0y$ for
some $w_0\in W$ and we have
\ben
\omega_{0,2}(x,y) = \lim_{y'\to y} \left(
\sum_{w\in W} (\gamma|w\gamma)B(x,wy') -
(\gamma|w_0\gamma)\frac{d\lambda(x)
  d\lambda(y')}{(\lambda(x)-\lambda(y'))^2} 
\right),
\een
where $y'$ is sufficiently close to $y$. Note that the set of poles of
$\omega_{0,2}(x,y)$ is the following set of points in $V_s\times V_s$:
\ben
\bigcup_{j=1}^N \ \lambda^{-1}(u_j)\times \lambda^{-1}(u_j).
\een

\subsection{Proof of Theorem \ref{t1}}

We are going to prove that the global recursion reduces to the local
one. To begin with let us simplify the kernel of the global
recursion. Put 
\ben
S_{y_1,y_2}(x) = \int_{y_2}^{y_1}B(x,y'),\quad x,y_1,y_2\in V_s.
\een
This is the unique form on $V_s$ with vanishing
$\mathcal{A}_i$-periods, with poles of order 1 at $y_1$ and
$y_2$ with residues respectively $+1$ and $-1$. The kernel of the local
recursion has the following symmetry
\ben
K_{\alpha,\beta} (s,x,wy) = K_{\alpha,w^{-1}\beta}(s,x,y).
\een
In particular, using this symmetry when $w=\theta_{j,a}$ and
Proposition \ref{K-Bergman} we get
\beq\label{K-S}
K_{\alpha,\beta_{j,a}}(s,x,y) = -\frac{1}{2} \sum_{w\in W}
(\alpha|w\beta_{j,a}) S_{w\theta_{j,a}(y),wy}(x).
\eeq
Furthermore, we have 
\ben
f_\gamma(y)-f_\gamma(\theta_{j,a}(y)) = \langle \gamma,
y-s_{\beta_{j,a}}(y)\rangle = \langle\beta_{j,a},y \rangle \,
(\gamma|\beta_{j,a}) .
\een
For fixed $j$, let us fix the local coordinate
$y_a:=\widetilde{\Phi}(C_a,s-\lambda\mathbf{1})$ 
near the ramification point $y_{j,a}$ ($C_a$ is a path
along which $\beta_{j,a}$ vanishes over $\lambda=u_j$). 
There is a unique element $w_a\in W$, s.t.,
$\beta_{j,a}=w_a\beta_{j,1}$ and $y_{j,a}=w_a y_{j,1}$ (recall that we
chose $\beta_{j,a}\in R_+$). We express the residue at a given
ramification point $y_{j,a}$ in terms of the residue at $y_{j,1}$. Let
us denote for brevity $\beta_j:=\beta_{j,1}$,
$\theta_j:=\theta_{j,1}$, and $y:=y_1$.  The
contribution on the RHS of the global recursion corresponding to the
$j$th term in the outer sum is
\ben
\frac{1}{2}\sum_{a=1}^{|W|/2} \operatorname{Res}_{y_a=y_{j,a}}
\frac{S_{\theta_{j,a}(y_a),y_a}(x_0)}{ (\gamma|\beta_{j,a}) \langle \beta_{j,a},y_a\rangle\,
  d\lambda}\times 
\left( \omega_{g-1,n+2}(s;y_a,\theta_{j,a}(y_a),x_1,\dots,x_n) +\cdots\right),
\een
where the omitted term differs from the corresponding term on the RHS
of the global recursion via the substitution $y\mapsto y_a$. After
changing the variables $y_a=w_a y$, the residue turns into a residue
at $y_{j,1}$, i.e., 
\ben
\frac{1}{2}\sum_{a=1}^{|W|/2} \operatorname{Res}_{y=y_{j,1}}
\frac{S_{w_a\theta_j(y),w_a y}(x_0)}{ (\gamma|\beta_{j,a}) \langle \beta_j,y \rangle\,
  d\lambda}\times 
\left( \omega_{g-1,n+2}(s;w_a y ,w_a \theta_{j}(y),x_1,\dots,x_n) +\cdots\right).
\een 
The term in the bracket  is by definition
\beq\label{rec-corr}
\omega^{w_a^{-1}\gamma,s_{\beta_j}w_a^{-1}\gamma,\gamma,\dots,\gamma}
_{g-1,n+2}(s;\lambda ,\lambda,\lambda_1,\dots,\lambda_n) +\sum
\omega_{g',n'+1}^{w_a^{-1}\gamma,\gamma_{I'} }(\lambda,\lambda_{I'})\, 
\omega_{g'',n''+1}^{s_{\beta_j}w_a^{-1}\gamma,\gamma_{I''} }(\lambda,\lambda_{I''})
\eeq
where $\lambda$ and $\lambda_1,\dots,\lambda_n$ are the projections of
$y$ and $x_1,\dots,x_n$ on the base of the branched covering \eqref{cover}.
We have the following decomposition
\ben
w_a^{-1}\gamma=\gamma'+ (w_a^{-1}\gamma|\beta_j)\beta_j/2= \gamma'+(\gamma|\beta_{j,a}) \beta_j/2,\quad
s_{\beta_j}(w_a^{-1}\gamma) = \gamma'-(\gamma|\beta_{j,a}) \beta_j/2,
\een
where $\gamma'$ is a cycle invariant with respect to the local monodromy
around $\lambda=u_j$. The period vectors $\phi_{\gamma'}(s,\lambda;z)$ are
analytic near $\lambda=u_j$, so up to terms that are
analytic at $y=y_{j,1}$ we get that \eqref{rec-corr} coincides with 
\ben
\frac{1}{4}\,(\gamma|\beta_{j,a})^2\,
\Big(\omega^{\beta_j,-\beta_j,\gamma,\dots,\gamma}
_{g-1,n+2}(s;\lambda ,\lambda,\lambda_1,\dots,\lambda_n) +\sum
\omega_{g',n'+1}^{\beta_j,\gamma_{I'} }(\lambda,\lambda_{I'})\, 
\omega_{g'',n''+1}^{-\beta_j,\gamma_{I''} }(\lambda,\lambda_{I''})
\Big).
\een
Note that $(1/2)\operatorname{Res}_{y=y_{j,1}}
=\operatorname{Res}_{\lambda=u_j}$. To finish the proof we just need to compute the sum 
\ben
\sum_{a=1}^{|W|/2}
(\gamma|\beta_{j,a}) S_{w_a\theta_j(y),w_a y}(x_0) = \frac{1}{2}\sum_{w\in W}
(\gamma|w\beta_j) S_{w\theta_j(y),w y}(x_0)=-K_{\gamma,\beta_j}(s,x_0,y).
\een
It remains only to recall Proposition \ref{phase-form} 
\ben
K_{\gamma,\beta_j}(s,x_0,y) = \Omega(\mathbf{\phi}^\gamma_+(s,\lambda_0;z),\mathbf{f}^{\beta_j}_-(s,\lambda;z))
\een
and to recall that by definition
\ben
\langle \beta_j, y\rangle = (I^{(-1)}_{\beta_j}(s,\lambda),1).\qed
\een

\section{The free energies and the primary potentials in genus 0 and 1}

The main goal of this section is to compute the genus-0 and genus-1
free energies. However, let us first prove that the spectral curve is non-singular.

\subsection{Smoothness of the spectral curve}
The spectral curve is a branched covering \eqref{cover} of a smooth curves, so the
only singularities could be at the ramification points. The
ramification points of index 2 are easy to analyze, because locally the
covering  near such a point is equivalent to a covering defined by the
period map of an $A_1$-singularity. Therefore, we can reduce the proof
of the general case to the case of an $A_1$-singularity. The latter
case is straightforward, so we omit the details. It is more
interesting to prove the regularity at the ramification points $\lambda^{-1}(\infty)$. 

Let us first recall several properties of the Coxeter
transformations. Given a Coxeter transformation $\sigma$, all other
Coxeter transformations have the form $w\sigma w^{-1}$, $w\in W$ and
the set of all Coxeter transformations consist of $|W|/h$
elements. Note that the number of ramification points above
$\lambda=\infty$ is also $|W|/h$. By definition, the ramification points are
the solutions of the following equations in $\mathbb{P}^{N-1}$: 
\ben
t_a(X_1,\dots,X_N) = 0,\quad 1\leq a\leq N-1.
\een
We assign a ramification point $\xi=[\xi_1,\dots,\xi_N]$ to each
Coxeter transformation $\sigma$, by letting $\sum_{i=1}^N
\xi_i\omega_i\in \mathfrak{h}$ be an eigenvector of $\sigma$ with
eigenvalue $\eta:=e^{2\pi\sqrt{-1}/h}$. Recall the so called Coleman
lemma \cite{Co}: $\langle \alpha,\xi\rangle \neq 0$ 
for all $\alpha\in R$, i.e., each eigenvector $\xi$ is inside some
Weyl chamber. 
\begin{proposition}
The map that associates a ramification point to a Coxeter
transformation is a bijection.
\end{proposition}
\proof
Let us assume that $\sigma_1\xi=\eta \xi=\sigma_2 \xi$
for two Coxeter transformations $\sigma_1$ and $\sigma_2$. Since
the Weyl group acts faithfully on the set of Weyl chambers and 
$\sigma_1^{-1}\sigma_2$ fixes the Weyl chamber to which $\xi$ belongs,
we must have $\sigma_1=\sigma_2$. Since both sets have the same number
of elements, the map must be onto.
\qed

Assume now that $\xi=[0,\xi_1,\dots,\xi_N]\in V_s$ is a ramification
point. We may assume that $\xi_N=1$, so the ramification point is
in the affine chart $U_N:=\{X_N\neq 0\}\subset \mathbb{P}^N$. Let
$u_i=X_i/X_N$, $0\leq i\leq N-1$ be the affine coordinates of
$U_N$. The equation of $V_s\cap U_N$ can be written as 
\ben
t_a(u_1,\dots,u_{N-1},1) = s_a u_0^{d_a},\quad 1\leq a\leq N-1.
\een  
Using the Jacobian criterion, we get that we have to prove that the determinant 
\ben
\operatorname{det}
\begin{bmatrix}
\frac{\partial t_1}{\partial x_1} & \cdots & \frac{\partial  t_1}{\partial x_{N-1}} \\
\vdots & & \vdots\\
\frac{\partial t_{N-1}}{\partial x_1} & \cdots & \frac{\partial
  t_{N-1}}{\partial x_{N-1}} 
\end{bmatrix}
\een
is non-zero at $x=(\xi_1,\dots,\xi_N)$. Let us look at the larger
determinant
\ben
\operatorname{det}
\begin{bmatrix}
\frac{\partial t_1}{\partial x_1} & \cdots & \frac{\partial t_1}{\partial x_{N-1}} & \frac{\partial t_1}{\partial x_{N}} \\
\vdots & & \vdots & \vdots\\
\frac{\partial t_{N-1}}{\partial x_1} & \cdots & \frac{\partial t_{N-1}}{\partial x_{N-1}}  & \frac{\partial
  t_{N-1}}{\partial x_{N}}  \\
\frac{\partial t_N}{\partial x_1} & \cdots & \frac{\partial
  t_N}{\partial x_{N-1}} & \frac{\partial t_N}{\partial x_{N}}
\end{bmatrix}
 =
\prod_{\alpha\in R_+} \langle \alpha, x\rangle.
\een
Since the invariant polynomials are weighted homogeneous,
we have
\ben
\sum_{i=1}^N x_i\frac{\partial t_a}{\partial x_i} = d_a t_a.
\een
Therefore, when we evaluate the bigger determinant at
$x=(\xi_1,\dots,\xi_N)$ we may replace the last column by  $(0,\dots,0,
d_N t_N(\xi))^t$, where we used that $t_a(\xi)=0$ for $1\leq a\leq
N-1$. Again, the Coleman's lemma implies that the big determinant is
non-zero, so both $t_N(\xi)$ and the determinant that we are
interested in must be non-zero. 

Finally, let us point out that our argument proves that the
ramification points $\lambda^{-1}(\infty)$ are smooth for all $s\in
\mathbb{C}^{N-1}$.

\subsection{Genus-0 free energy}

The genus-0 free energy is defined through the meromorphic
differential 
\ben
f_\gamma(x) d\lambda = (I^{(-1)}_\gamma(s,\lambda),1) d\lambda = d_\lambda
(I^{(-2)}_\gamma(s,\lambda),1), 
\een
where $x=\widetilde{\Phi}(C,s-\lambda\mathbf{1})$ and $C$ is the path
to the reference point that determines the value of the period. 
The poles of this differential are only at the ramification points
$\{x_a\}_{a=1}^{|W|/h}:=\lambda^{-1}(\infty)$ and the integrals along any closed path in
$V_s$ is $0$, so the definition from \cite{EO} takes the form
\ben
\underline{F}^{(0)} = \frac{1}{2}\sum_{a=1}^{|W|/h}
\operatorname{Res}_{x=x_a} V_a(x) f_\gamma(x) d\lambda(x),
\een
where
\ben
V_a(x)=\operatorname{Res}_{y=x_a} \log(1-\zeta(x)/\zeta(y))\,  f_{\gamma}(y) \, d\lambda(y),
\een
where $\zeta:U_a-\{x_a\}\to \mathbb{C}$ is a local coordinate in a neighborhood
$U_a$ of $x_a$, s.t., $\lambda(y) = \zeta(y)^h$ $\forall y\in U_a$.

We have $\mathbf{f}^\alpha(s,\lambda;z)=S_s(z)
\mathbf{f}^\alpha(0,\lambda;z)$, where $S_s=1+S_1 z^{-1}+\cdots$ is a
fundamental solution for the Dubrovin's connection 
\ben
z\partial_{t_a} S_t(z) = v_a\bullet S_t(z),\quad S_0(z)=1. 
\een
Let us denote by $\sigma:\mathfrak{h}\to \mathfrak{h}$ the Coxeter
transformation corresponding to the monodromy along a big loop around
the discriminant (in counterclockwise direction), then using the
homogeneity of $v_i$, we get 
\ben
(I^{(0)}_\alpha(0,\lambda),v_i) = \lambda^{-m_i/h}\langle H_i,\alpha\rangle,
\een 
where $m_i=d_i-1$ are the Coxeter exponents and $H_i\in \mathfrak{h}$ is an
eigenvector of $\sigma$ with eigenvalue $\eta^{m_i}$. Note that $m_1=1$, $m_N=h-1$,
$m_i+m_{N+1-i}=h$. Using Saito's  formula \eqref{Saito-formula} we
get that the eigenbasis $\{H_i\}_{i=1}^N$ satisfy 
\ben
(H_i|H_j) = \delta_{i+j,N+1},\quad 1\leq i,j\leq N.
\een
It follows that we can express the free genus-0 energy in terms of the
eigenbasis and the matrices $S_k$. After a direct computation we get
the following formula for $(I^{(-2)}_\gamma(s,\lambda),1)$
\ben
\sum_{i=1}^N
\frac{\lambda^{-m_i/h+2}}{(-m_i/h+1)(-m_i/h+2)}\langle\gamma,H_i\rangle
(v^i,1) -
\sum_{i=1}^N
\frac{\lambda^{-m_i/h+1}}{-m_i/h+1}\langle\gamma,H_i\rangle
(S_1v^i,1) +\\
+\sum_{k=0}^\infty \sum_{i=1}^N (m_i/h)(m_i/h+1)\cdots
(m_i/h+k-1)\lambda^{-m_i/h-k}
\langle \gamma, H_i\rangle (S_{k+2}v^i,1),
\een 
where $\{v^i\}$ is a basis of $H$ dual to $\{v_i\}$ with respect to
the residue pairing.

Let us assume first that $x_a$ is the ramification point corresponding
to the classical monodromy, then 
\ben
V_a(x) = \zeta^{h+1}\frac{h^2}{h+1}\langle \gamma, H_N\rangle -
\sum_{i=1}^N\frac{\zeta^{m_i}}{m_i}\langle \gamma, H_{N+1-i}\rangle (S_1v^{N+1-i},1),
\een
where $\zeta=\zeta(x)=\lambda(x)^{1/h}$ is the local coordinate near
$x_a$. From this formula we get that
\ben
\operatorname{Res}_{x=x_a} V_a(x) f_\gamma(x)d\lambda(x) = 
-\operatorname{Res}_{\zeta=\infty} (I^{(-2)}_\gamma(s,\zeta^h),1)\,
d_\zeta V_a(x) 
\een
is 
\beq\label{genus0-H}
h\Big( \langle \gamma, H_1\rangle \, \langle \gamma, H_N\rangle (S_3 1,1) - 
\sum_{i=1}^N \langle \gamma, H_i\rangle \, \langle \gamma,
H_{N+1-i}\rangle (S_2v^i,1)(S_1 v^{N+1-i},1)
\Big).
\eeq
The above formula can be simplified as follows. Note that 
\ben
h \langle \gamma, H_i\rangle \, 
H_{N+1-i} = \sum_{k=1}^h \eta^{m_i k}\sigma^k\gamma,
\een
so \eqref{genus0-H} takes the form
\beq\label{genus0-Cox}
\sum_{k=1}^h \eta^{ k}\langle \sigma^k\gamma|\gamma) (S_31,1) - 
\sum_{k=1}^h\sum_{i=1}^N \eta^{m_i k}\langle \sigma^k\gamma|\gamma)
(S_2v^i,1)(S_1 v^{N+1-i},1).
\eeq
The expression \eqref{genus0-Cox} is the contribution to
$2\underline{F}^{(0)}$ coming from the residue at the ramification
point $x_a$. Note that after adding the remaining contributions we get
that $2\underline{F}^{(0)}$ is the sum of \eqref{genus0-Cox} over all
Coxeter transformations $\sigma$.
\begin{lemma}
The following identity holds
\ben
\sum_{k=1}^h \eta^{m_i k}\sum_\sigma \sigma^k = |W|/N,
\een
where the sum is over all Coxeter transformations $\sigma$.
\end{lemma}
\proof
The operator $\sum_\sigma \sigma^k$ commutes with the action of $W$,
so By Schur's lemma, it must act by some constant $c_k$. After taking trace we get 
\ben
c_k N = \operatorname{Tr}(\sigma^k) |W|/h,
\een
where we used that there are $|W|/h$ Coxeter transformations and that
the trace of $\sigma^k$ is the same for all Coxeter
transformations. On the other hand, 
\ben
\sum_{k=1}^h \eta^{m_i k} \operatorname{Tr}(\sigma^k) =
\sum_{j=1}^N\sum_{k=1}^h \eta^{(m_i-m_j)k} = h.\qed
\een
Applying the above Lemma and using that $S_t(z)S_t(-z)^T=1$, we get 
\ben
\underline{F}^{(0)} = \frac{1}{2} ( (S_3-S_2 S_1)1,1)\, |W|(\gamma|\gamma)/N.
\een
Using that $S_t(z)$ is a solution for the Dubrovin's connection, it is
easy to verify that 
\ben
F^{(0)} = \frac{1}{2}((S_2 S_1-S_3)1,1)
\een
is a potential of the Frobenius structure, so up to quadratic terms in
$t$ we have 
\ben
\underline{F}^{(0)}(t) = -(\gamma|\gamma)\, \frac{|W|}{N} F^{(0)}(t). 
\een
\subsection{Genus-1 free energy}
Let us denote by $x_{j,a}$ ( $1\leq j\leq N$, $1\leq a\leq |W|/2$) the
double ramification points and by $u_j:=\lambda(x_{j,a})$ the
corresponding branching points.  The genus-1 free energy is by definition
\ben
\underline{F}^{(1)}(s) = -\frac{1}{2}\log \tau_B(s)
-\frac{1}{24}\sum_{i=1}^N\sum_{a=1}^{|W|/2} \log f_\gamma'(x_{j,a}),
\een
where $f_\gamma'$ is the derivative with respect to the local 
parameter $\sqrt{\lambda(x)-u_j}$ near $x=x_{j,a}$ and $\tau_B$ is the
Bergman tau-function of $V_s$. 

Recall that the period mapping has the following Laurent series
expansion near $\lambda=u_j$ (see \cite{G3, M1}):
\ben
I^{(0)}_{\beta_j}(s,\lambda) = \pm
\frac{2}{\sqrt{2(\lambda-u_j)\Delta_j}}
\Big( du_j+\cdots\Big),
\een
where the dots represents higher order terms in $\lambda-u_j$,
$\Delta_j:=(du_j,du_j)$, and $\beta_j$ is a cycle vanishing over
$\lambda=u_j$. The points $x=\widetilde{\Phi}(C,s-\lambda\mathbf{1})$
in a neighborhood of $x_{j,a}$ correspond to $\lambda$ in a
neighborhood of $u_j$, so we get
\ben
f_\gamma'(x_{j,a}) = \lim_{\lambda\to u_j} \frac{d_\lambda
  (I^{(-1)}_\gamma(s,\lambda),1)  }{d_\lambda \sqrt{\lambda-u_j} } = 
\lim_{\lambda\to u_j}  2\sqrt{\lambda-u_j}\, 
  (I^{(0)}_\gamma(s,\lambda),1) .
\een
Decomposing the cycle $\gamma = \gamma'+(\gamma|\beta_j)\beta_j/2$
into invariant and anti-invariant parts with respect to the local
monodromy we get 
\ben
f_\gamma'(x_{j,a})=\lim_{\lambda\to u_j}  (\gamma|\beta_j)
\sqrt{\lambda-u_j} \,  \frac{\pm 2}{\sqrt{2(\lambda-u_j)\Delta_j}} =
\pm\sqrt{2}\, (\gamma|\beta_j)\,\Delta_j^{-1/2}. 
\een
\subsubsection{The Bergman $\tau$-function}
To define the Bergman tau-function we have to think of the pair
$(V_s,\lambda)$ as a point in an appropriate moduli space
$\mathcal{M}$ of Hurwitz covers of $\mathbb{P}^1$ whose genus and
ramification profile is the same as of 
$V_s$. The critical values $u_{j,a}=\lambda(x_{j,a})$ provide local
coordinates on $\mathcal{M}$ and the differential of $\tau_B$ at
$(V_s,\lambda)$ is defined via  
\ben
d\log \tau_B = \sum_{j=1}^N \sum_{a=1}^{|W|/2} du_{j,a} \,
\operatorname{Res}_{x=x_{j,a}} \frac{B(x,\theta_{j,a}(x))}{d\lambda},
\een
where $B$ is the Bergman kernel of $V_s$.

Let $u=(\lambda(x)-u_j)^{1/2}$ and $v=(\lambda(y)-u_j)^{1/2}$ be the
local coordinates of two points $x,y\in V_s$ near $x_{j,a}$. The
Bergman kernel has the form
\ben
B(x,y) = \frac{du dv}{(u-v)^2} + f_{j,a}(u,v)dudv,
\een
where $f_{j,a}\in \mathbb{C}\{u,v\}$ is a convergent power series in $u$ and $v$. If
$y=\theta_{j,a}(x)$, then $v=-u$ and we get that 
\ben
\operatorname{Res}_{x=x_{j,a}} \frac{B(x,\theta_{j,a}(x))}{d\lambda} =
-\frac{1}{2}f_{j,a}(0,0). 
\een
\subsubsection{The Bergman $\tau$-function and the $R$-matrix}
Following the notation in \cite{M1} we recall the following formula
for the correlator
\ben
\omega_{0,2}^{\beta_j,\beta_j}(s;\lambda,\lambda)=P_0^{jj}(s,\lambda)d\lambda\cdot
d\lambda,
\een
where $\beta_j$ is the vanishing cycle vanishing over $\lambda=u_j$:
\ben
P_0^{jj}(s,\lambda) = \frac{1}{4} (\lambda-u_j)^{-2}
+ 2(e_j,V_{00}(s)e_j) (\lambda-u_j)^{-1},
\een
where $e_j=du_j/\sqrt{\Delta_j}\in T_s^*B\cong H$ and $V_{k\ell}(s)$ are linear
operators of $H$ defined via Givental's $R$-matrix
$\mathcal{R}(s,z)=1+R_1(s)z+R_2(s)z^2+\cdots$ 
\ben
\sum_{k,\ell=0}^\infty V_{k\ell}(s) w^kz^\ell =
\frac{1-{}^T\mathcal{R}(s,-w)\mathcal{R}(s,-z)}{z+w}. 
\een
Note that $V_{00}=R_1$ and since $\{e_j\}$ is an orthonormal basis of
$H$, we get that $(e_j,V_{00}(s)e_j)=R_1^{jj}(s)$ is the $j$th diagonal
entry of $R_1(s)$. 

Recalling the definition of the 2-point genus-0 correlators we get that
\beq\label{residue}
\operatorname{Res}_{x=x_{j,a} }
\frac{1}{d\lambda}\left. \Big(
  \Omega(\phi_+^\gamma(s,\lambda;z),\phi_+^\gamma(s,\mu;z)
  -(\gamma|\gamma)\frac{d\lambda\cdot d\mu}{(\lambda-\mu)^2} \Big)
  \right|_{\mu=\lambda} 
\eeq
is 
\ben
&&
 \operatorname{Res}_{x=x_{j,a} }
 P^{(0)}_{\gamma,\gamma}(s,\lambda)d\lambda =
 \operatorname{Res}_{\lambda=u_j} 
\Big( 
P^{(0)}_{\gamma,\gamma}(s,\lambda)
+P^{(0)}_{\theta_{j,a}\gamma,\theta_{j,a}\gamma}(s,\lambda) 
\Big)d\lambda
=\\
&&
\operatorname{Res}_{\lambda=u_j}
\Big(
2P^{(0)}_{\gamma',\gamma'}(s,\lambda)+
\frac{1}{2} (\gamma|\beta_j)^2 \, P_0^{jj}(s,\lambda)\Big)d\lambda = 
 (\gamma|\beta_j)^2 R_1^{jj}(s),
\een
where $x_{j,a}$ is the ramification point corresponding to the
reference path that defines the residue \eqref{residue} and $\gamma'$ is
the invariant part of $\gamma$ with respect to the local monodromy
around $\lambda=u_j$. Note that $P_{\gamma',\gamma'}^{(0)}$ is
holomorphic near $\lambda=u_j$, so it does not contribute to the
residue. 

On the other hand we can compute the residue \eqref{residue} in terms
of the Bergman kernel. Recalling Proposition \ref{phase-form} and
Proposition \ref{K-Bergman} we transofrm the residue \eqref{residue}  into
\ben
\operatorname{Res}_{x=x_{j,a} }\left.\left(
\sum_{w\in W} (\gamma |w\gamma) \, \frac{B(x,wy)}{d\lambda(x)} -
(\gamma|\gamma) \frac{d\lambda(y)}{(\lambda(x)-\lambda(y))^2}
\right)\right|_{y=x}.
\een
The only terms in the above sum that contribute to the residue are the ones for which
$w=1$ or $w=\theta_{j,a}$. Using again the local coordinates
$u=\sqrt{\lambda(x)-u_j}$ and $v=\sqrt{\lambda(y)-u_j}$ we get that
the term with $w=1$ and the term outside of the sum add up to 
\ben
(\gamma|\gamma)\Big( \frac{1}{8} u^{-3} + \frac{1}{2}f_{j,a}(u,u) u^{-1}\Big) du
\een  
and the contribution to the residue is $(\gamma|\gamma)
f_{j,a}(0,0)/2$. The term with $w=\theta_{j,a}$ contributes to the
residue
\ben
-\frac{1}{2} (\gamma|\theta_{j,a}\gamma) f_{j,a}(0,0). 
\een
Comparing the two computations of the residue \eqref{residue}  we get
\ben
(\gamma|\beta_j)^2 R_1^{jj}(s)=
\frac{1}{2}(\gamma|\gamma-\theta_{j,a}(\gamma)) f_{j,a}(0,0)= 
\frac{1}{2}(\gamma|\beta_j)^2 f_{j,a}(0,0),
\een
i.e., $\frac{1}{2} f_{j,a}(0,0)=R_1^{jj}$. Therefore the differential
of the Bergman $\tau$-function is the following 
\ben
d\log\tau_B = -\frac{|W|}{2} \sum_{j=1}^N R_1^{jj}(s) du_j. 
\een
Finally, the genus-1 free energy becomes
\ben
\underline{F}^{(1)}(s) = \frac{|W|}{2} \left(
\frac{1}{2}\int \sum_{j=1}^N R_1^{jj}(s) du_j + \frac{1}{48}
\sum_{j=1}^N \log \Delta_j\right).
\een
It remains only to recall that the term in the brackets is the genus-1
primary potential of the Frobenius structure also known as the
$G$-function (see \cite{G2} for more details and references).

\medskip

\noindent 
{\bf Acknowledgments}. 
This work is supported in part by JSPS Grant-in-Aid 26800003.
We acknowledge the World Premiere International Research Center Initiative (WPI Initiative), Mext, Japan. 

\appendix

\section{The space of holomorphic 1-forms}\label{App1}

In this appendix we would like to give a description of the space of
holomorphic 1-forms and to prove that our
initial condition \eqref{w02} is independent of the choice of a
Torelli marking. Unfortunately, our argument does not work in the
exceptional cases. Of course, we can use Corollary
\ref{K-Bergman:cor}, but it would be nice to find a direct algebraic
proof.  We also give an amusing proof that the order of the
Weyl group is the product of the degrees of the invariant
polynomials. 

\subsection{The space of holomorphic 1-forms}
Using
the Riemann--Hurwitz formula we can compute the genus of $V_s$ as 
\ben
g=1+ d(V) |W|/(2h),
\een
where
\ben 
d(V) = Nh/2-h-1=d_1+d_2+\cdots+d_{N-1}-N-1.
\een
The genus can be computed also using that $V_s$ is a complete
intersection, which allows us to compute the canonical bundle via the
adjunction formula (see \cite{Di1})
\ben
g=1+d_2\cdots
d_{N-1}(1+d_2+\cdots+d_{N-1}-N) . 
\een
In particular, we get a uniform proof that $|W|=d_1d_2\dots d_N$. 

The space of holomorphic 1-forms on $V_s$ can be described as follows
(see \cite{Di2}).  In the affine chart $X_0\neq 0$ and
the affine coordinates $x_i=X_i/X_0$ it is easy to see that 
\beq\label{hol-forms}
\phi(x_1,\dots,x_N) \frac{dx_1\wedge \cdots\wedge dx_N}{dt_1\wedge \cdots\wedge dt_{N-1}}
\eeq
extends to a holomorphic form on $V_s$ if and only if $\phi\in
\mathbb{C}[x_1,\dots, x_N]$ is a polynomial of degree at most
$d(V)$. Note that the form \eqref{hol-forms} is identically 0 on $V_s$
if and only if $h\in (t_1(x)-s_1,\dots,t_{N-1}(x)-s_{N-1})$. It
remains only to check that the number of elements in the ring
\ben
\mathbb{C}[x_1,\dots, x_N]/(t_1(x)-s_1,\dots,t_{N-1}(x)-s_{N-1}).
\een
of degree at most $d(V)$ is $g$.
\begin{proposition}\label{av:identity}
If $\phi\in \operatorname{Sym}(\mathfrak{h}^*) $ is a polynomial of degree at
most $d(V)$ and $\gamma\in \mathfrak{h}^*$ is a linear function, then
\ben
\sum_{w\in W} \operatorname{det}(w)\, (w\gamma \otimes w^{-1} \phi) = 0.
\een
\end{proposition}
\proof
Our argument works in the $A$ and $D$ cases only. The identities in
the exceptional cases, can be verified with a computer. 

The LHS will be viewed as a function $f$ on $\mathfrak{h}\times
\mathfrak{h}$
\ben
f(x,y) = \sum_{w\in W} \operatorname{det}(w) \langle
w\gamma,x\rangle\, \phi(w y),\quad (x,y)\in \mathfrak{h}\times \mathfrak{h}. 
\een 
Since the function depends linearly on $\gamma$, it is enough to prove
the identity for a set of $\gamma$'s that form a basis of
$\mathfrak{h}^*$. Similarly, we may assume that $\phi$ is a monomial
in $y$. 

Let us take $\gamma$ to be a fundamental weight
corresponding to a node of the Dynkin diagram, s.t., if we remove that
node, then we get a Dynkin diagram of the same type but with rank one
less. Note that the number of positive roots orthogonal to $\gamma$ is
$\frac{1}{2} (N-1)N $ for $A_N$ and $(N-2)(N-1)$ for $D_N$. In both
cases, the number is greater than 
\ben
d(V)=
\begin{cases}
\frac{1}{2}N(N-1) -2 & \mbox{ for } A_N\\
(N-2)(N-1)-1 & \mbox{ for } D_N.
\end{cases}
\een 
In particular, the polynomial
\ben
\Delta_\gamma(y)=\prod_{\alpha\in R_+: (\alpha|\gamma)=0} \langle
\alpha ,y\rangle
\een
has degree at least $d(V)+1$. The zero locus of $\Delta_\gamma$ is
contained in the zero locus of $f$: if $\langle \alpha ,y_0\rangle =0
$, then in the definition of $f$ let us shift the summation by
replacing $w\mapsto ws_\alpha$, we get $f(x,y_0)=-f(x,y_0)$. The ideal
generated by $\Delta_\gamma$ is a radical ideal, so using Hilbert's
Nullstellensatz, we get that $f(x,y)=g(x,y)\Delta_\gamma(y)$ for some
polynomial $g$. If we assume that $f\neq 0$, then we get a
contradiction by comparing the degrees of the monomials in $y$ on both
sides: on the left they all have degree $\operatorname{deg}(\phi)\leq
d(V)$, while on the right, they all have degree at least
$\operatorname{deg}(\Delta_\gamma)>d(V)$. 

To finish the proof, we just need to use that the above argument
applies to the entire orbit $W\gamma$ and that this orbit contains a
basis of $\mathfrak{h}^*$. 
\qed.  

The above proposition implies that our initial condition is
independent of the Torelli marking. Indeed, changing the Torelli
marking will modify the Bergman kernel via a quadratic expression of
holomorphic differentials on $V_s$. Using the explicit description of
the holomorphic differentials from above we see that if we replace $B(x,y)$
in \eqref{w02} by a product $\theta_1(x)\theta_2(y)$ of holomorphic
differentials, then  we get precisely the identity in Proposition \ref{av:identity}.

\end{document}